\documentclass{aip-cp}

\usepackage[numbers]{natbib}
\usepackage{rotating}
\usepackage{graphicx}
\usepackage{verbatim}

\def\ds{\displaystyle}
\def\eps{\varepsilon}
\def\A{\mathbf{A}}
\def\CC{\mathbb{C}}
\def\II{\mathbf{I}}
\def\O{\mathcal{O}}
\def\M{\mathbf{M}}
\def\N{\mathbf{N}}

\def\PP{\mathbf{P}}
\def\PPhi{\mathbf{\Phi}}

\newcommand{\eqref}[1]{{(\ref{#1})}}
\newtheorem{theorem}{Theorem}  

\begin{document}

% Title portion
\title{High-order WKB-based Method For The 1D Stationary Schr\"odinger Equation In The Semi-classical Limit}

\author[aff1]{Anton Arnold\corref{cor1}}
\eaddress[url]{https://www.asc.tuwien.ac.at/arnold/}
\author[aff1]{Jannis K\"orner}
\eaddress{jannis.koerner@tuwien.ac.at}

\affil[aff1]{Inst.\ f.\ Analysis u.\ Scientific Computing, Technische Universit\"at Wien,  Wiedner Hauptstr. 8, A-1040 Wien, Austria.}
\corresp[cor1]{Corresponding author: anton.arnold@tuwien.ac.at}

\maketitle

\begin{abstract}
We consider initial value problems for $\varepsilon^2\,\varphi''+a(x)\,\varphi=0$ in the highly oscillatory regime, i.e., with $a(x)>0$ and $0<\varepsilon\ll 1$. We discuss their efficient numerical integration on coarse grids, but still yielding accurate solutions. The $\O(h^2)$ one-step %marching 
method from \cite{ABN} is based on an analytic WKB-preprocessing of the equation. Here we extend this method to $\O(h^3)$ accuracy. 
\end{abstract}

% Head 1
\section{INTRODUCTION}

This paper is concerned with efficient numerical methods for highly oscillatory ordinary differential equations (ODEs) of the form
\begin{equation} \label{SchIVP}
\displaystyle \varepsilon^2 \varphi''(x) + a(x) \varphi(x)= 0 \,, \quad x \in (0,1)\,;\qquad
\ds \varphi(0)=\varphi_0\in\CC\,,\quad
\ds \varepsilon \varphi'(0)=\varphi_1\in\CC \,.
\end{equation}
Here, $0<\eps \ll 1$ is a small parameter and $a(x)\ge a_0>0$ a sufficiently smooth function, such that \eqref{SchIVP} does not include a turning point. For extensions with a turning point, i.e.\ a sign change of $a(x)$, we refer to \cite{AN, AD}.
Such problems have applications, e.g.\ in quantum transport \cite{BAP06, Cla_Naoufel}, mechanical systems (see references in \cite{lub}), and cosmology \cite{AHLH20}.

For $\eps\ll 1$, solutions to \eqref{SchIVP} are highly oscillatory, and hence standard ODE-solvers become inefficient since they need to resolve each oscillation by choosing $h=\O(\eps)$. In \cite{lub}, an $\eps$-uniform scheme with $\O(h^2)$ accuracy for large step sizes up to $h =\O(\sqrt\eps)$ was constructed, see also \S XIV of \cite{HLW} and references therein. The $\O(h^2)$-scheme of \cite{ABN} is based on a (w.r.t.\ $\eps$) second order WKB-approximation of \eqref{SchIVP} and makes the method even \emph{asymptotically correct}, i.e.\ the error decreases with $\eps$ even on a coarse spatial grid, if the phase function can be obtained analytically or with spectral accuracy \citep{AKU}. Here we present an $\O(h^3)$ extension of the latter method; for its detailed analysis we refer to \cite{ADK}.

%%%%%%%%%%%%%%%%%%%%%%%%%%%%%%%%%%%%%%%%%%%%%%%%%%%%%%%%%%%%%%

\section{WKB-TRANSFORMATION AS ANALYTIC PREPROCESSING}

The essence of this numerical method is to transform the highly oscillatory problem \eqref{SchIVP} into a much ``smoother'' problem by eliminating the dominant oscillation frequency. Following \cite{ABN} we first introduce the vector function 
$U(x):=\Big(a^{1/4} \varphi(x)\,,\,\ds \frac{\varepsilon (a^{1/4}\varphi)'(x)}{\sqrt{a(x)}}\Big)^\top$. Then we set 
$Z(x) := e^{-{i\over \eps} \PPhi^\eps(x)} \PP\,U(x)$ with the matrices
$$\PP := {1\over \sqrt{2}}
\left(
\begin{array}{cc}
i&1\\
1&i
\end{array}
\right)\;; \quad \PP^{-1} = {1\over \sqrt{2}}
\left(
\begin{array}{cc}
-i&1\\
1&-i
\end{array}
\right)\,,
$$
\begin{equation} \label{phase}
\mathbf{\Phi}^\eps (x) := 
%\int_0^x D^\eps(y)\, dy = 
\mbox{diag}(\phi(x), -\phi(x))\,;\qquad
\phi(x):=\int_0^x \left( \sqrt{a(\tau)} - \eps^2 b(\tau)\right) \,d\tau\,; \qquad 
b(x) := -\frac{1}{2a(x)^{1/4}} \big(a(x)^{-1/4}\big)''\,.
\end{equation}
We remark that the (real valued) phase function $\phi$ is precisely the phase in the (w.r.t.\ $\eps$) second order WKB-approximation of \eqref{SchIVP} (cf.\  \cite{ABN,LL85}).
Then, $Z$ satisfies the ODE initial value problem (IVP)
\begin{equation} \label{EQZ}
\ds Z' = \eps \N^\eps(x) Z\,,\quad x\in(0,1);\qquad \ds Z(0)=Z_I=\PP\,U_I\,; \qquad U_I= U(0)\,.
\end{equation} 
$\N^\eps$ is an off-diagonal matrix with the entries 
$N^\eps_{1,2} (x)= b(x) 
e^{-\frac{2i}{ \eps} \phi(x)},\,
N^\eps_{2,1} (x)= b(x) e^{\frac{2i}{ \eps} \phi(x)}$.
While the ODE \eqref{EQZ} is still oscillatory, in fact with doubled frequency, $Z$ is ``smoother'' than $\varphi$ and $U$, as its oscillation amplitude is reduced to $\O(\eps^2$), cf.\ \cite{ABN}. 
After numerically solving the ODE \eqref{EQZ}, the original solution is recovered by
$U(x)=\PP^{-1} e^{{i\over \eps} \PPhi^\eps(x)} Z(x)\,.$

%%%%%%%%%%%%%%%%%%%%%%%%%%%%%%%%%%%%%%%%%%%%%%%%%%%%%%%%%%%%%%

\section{ASYMPTOTICALLY CORRECT NUMERICAL SCHEME}

To construct an asymptotically correct one-step scheme for the IVP \eqref{EQZ} on the uniform grid $x_n:=n\,h;\,n=0,...,N$ with the step size $h=1/N$, we consider first the truncated Picard iteration for \eqref{EQZ} (with $P=2$ in \cite{ABN}, and $P=3$ for the $O(h^3)$ method here):
$$
%	\begin{equation}\label{picard_limit}
		Z(\eta)\approx Z(\xi)+\sum_{p=1}^{P}\eps^{p}\mathbf{M}_{p}^{\eps}(\eta;\xi)\,Z(\xi)\,,
%	\end{equation}
$$	
where the matrices $\mathbf{M}_{p}^{\eps}$, $p=1,\,2,\,3$ are given by the iterated oscillatory integrals
$$
%	\begin{equation}
%		\mathbf{M}_{p}^{\eps}(\eta;\xi)&=\int_{\xi}^{\eta}\int_{\xi}^{y_{1}}\cdots\int_{\xi}^{y_{p-1}}\mathbf{N}^{\eps}(y_{1})\cdots\mathbf{N}^{\eps}(y_{p})\,dy_{p}\cdots\mathrm{d}y_{1}\Comma\nonumber\\
		\mathbf{M}_{p}^{\eps}(\eta;\xi)=\int_{\xi}^{\eta}\mathbf{N}^{\eps}(y)\mathbf{M}_{p-1}^{\eps}(y;\xi)\,dy\,,\quad \mathbf{M}_{0}^{\eps}=\II\,.
%	\end{equation}
$$
	
This is followed by a high order approximation of $\M_p^\eps$ (w.r.t.\ both small parameters $h$ and $\eps$) using the \emph{asymptotic method} for oscillatory integrals \cite{INO} and a shifted variant \cite{ABN}. We denote these approximation matrices by $\A_n^{p,P}\approx \eps^p\M_p^\eps(x_{n+1};x_n);\,p=1,...,P$. The two resulting numerical schemes, referred to as WKB2 (for $P=2$) and WKB3 (for $P=3$) have the structure: 
$$
%\begin{equation}\label{scheme}
	Z_{n+1}:=\left(\mathbf{I}+\sum_{p=1}^P \mathbf{A}^{p,P}_{n}\right)Z_{n}\,,\quad n=0,\dots,N-1\,.
%\end{equation}
$$
For the coefficients of $\A_n^{p,P}$ we have:
$$
%	\begin{equation}\label{Q1}
		\mathbf{A}_{n}^{1,P}:=\eps \left(
\begin{array}{cc}
			0 & \overline{Q_{1}^{P}(x_{n+1},x_n)} \\
			Q_{1}^{P}(x_{n+1},x_n) & 0 
		\end{array}\right)\,,\quad
		\mathbf{A}_{n}^{2,P}:=\eps^2 \left(
\begin{array}{cc}
			Q_{2}^{P}(x_{n+1},x_n) & 0 \\
			0 & \overline{Q_{2}^{P}(x_{n+1},x_n)} 
		\end{array}\right)\,,\quad P=2,\,3\,,	
%	\end{equation}
$$
with
\begin{eqnarray*}
	Q_{1}^{P}(x_{n+1},x_n) \!\!\!\!\!\!\!&&:=-\displaystyle	\sum_{p=1}^{P}(i\eps)^{p}\left(b_{p-1}(x_{n+1}) e^{\frac{2 i}{\eps}\phi(x_{n+1})}-b_{p-1}(x_{n}) e^{\frac{2 i}{\eps}\phi(x_{n})}\right)
- e^{\frac{2 i}{\eps}\phi(x_{n})}\displaystyle	\sum_{p=1}^{P}(i\eps)^{p+P}b_{p+P-1}(x_{n+1})\,
h_{p}\Big(\frac{2}{\eps}s_n\Big)\,,\\[3mm]
\ds  Q_{2}^{2}(x_{n+1},x_n) \!\!\!\!\!\!\!&&:=  \ds 
- i \eps (x_{n+1} -x_n) { b(x_{n+1})b_0(x_{n+1}) +b(x_{n}) b_0(x_{n}) \over 2} \\
&&  \ds - \eps^2 b_0(x_{n}) b_0(x_{n+1})\, 
 h_1\Big(-{2\over \eps}s_{n}\Big)
 +i  \eps^3 b_1(x_{n+1})  [b_0(x_{n})-b_0(x_{n+1})]\,  h_2\Big(-{2\over \eps}s_{n}\Big)\,,\\[3mm]
	Q_{2}^3(x_{n+1},x_n)\!\!\!\!\!\!\!&&:=-i\eps Q_{S}[bb_{0}](x_{n+1},x_n)\\
	&&-\eps^{2}\Big[b_{0}(x_n) e^{\frac{2 i}{\eps}\phi(x_n)}\left[b_{0}(y) e^{-\frac{2 i}{\eps}\phi(y)}\right]_{x_n}^{x_{n+1}}-Q_{S}[bb_{1}](x_{n+1},x_n)\Big]\\
	&&+ i\eps^{3}\big[b_{0}(x_n)b_{1}(x_{n+1})-b_{1}(x_n)b_{0}(x_{n+1})\big]\,
	h_{1}\Big(-\frac{2}{\eps}s_{n}\Big)\\
	&&+\eps^{4}\big[\left(b_{0}(x_n)+b_{0}(x_{n+1})\right)b_{2}(x_{n+1})-b_{1}(x_n)b_{1}(x_{n+1})-2b_{0}(x_{n+1})b_{3}(x_{n+1})s_{n}\big]\, 
	h_{2}\Big(-\frac{2}{\eps}s_{n}\Big)\\
	&&+ i\eps^{5}\big[\left(b_{0}(x_{n+1})-b_{0}(x_n)\right)b_{3}(x_{n+1})-\left(b_{1}(x_{n+1})-b_{1}(x_n)\right)b_{2}(x_{n+1})\big]\,
	h_{3}\Big(-\frac{2}{\eps}s_{n}\Big)\,,
\end{eqnarray*}
and the abbreviations
$$
    s_n:= \phi(x_{n+1})-\phi(x_n)\,;\quad
	b_{0}(x):=\frac{b(x)}{2\phi^{\prime}(x)}\,,\quad b_{p}(x):=\frac{b_{p-1}^{\prime}(x)}{2\phi^{\prime}(x)}\,; \qquad 
	h_{p}(x):= e^{i x}-\sum_{k=0}^{p-1}\frac{(i x)^{k}}{k!}\,, \quad p=1,\,2,\,3\,,
$$
$$
  	Q_{S}[f](\eta,\xi):=\frac{\eta-\xi}{6}\left(f(\xi)+4f\left(\frac{\xi+\eta}{2}\right)+f(\eta)\right)\,.
$$
%Since the off-diagonal elements of $\mathbf{A}_{n}^{3,3}$ are quite lengthy, we don't list them here, but refer to \cite{ADK}.
Finally we have 
$$
%	\begin{equation}\label{Q3}
		\mathbf{A}_{n}^{3,3}:=\eps^3 \left(
\begin{array}{cc}
			0 & \overline{Q_{3}^{3}(x_{n+1},x_n)} \\
			Q_{3}^{3}(x_{n+1},x_n) & 0 
		\end{array}\right)\,,	
%	\end{equation}
$$
with
\begin{eqnarray*}
	Q_{3}^{3}(x_{n+1},x_n) \!\!\!\!\!\!\!&&:
	=-\eps^{2} e^{\frac{2 i}{\eps}\phi(x_n)}\Bigg[\frac{x_{n+1}-x_n}{2}\left[c_{0}(x_{n+1})+b(x_n)b_{0}(x_n)b_{0}(x_{n+1})\right]\,h_{1}\Big(\frac{2}{\eps}s_{n}\Big)\Bigg] \\
		&&- i\eps^{3} e^{\frac{2 i}{\eps}\phi(x_n)}\Bigg[\frac{1}{2}\left[c_{1}(x_{n+1})(x_{n+1}-x_n)+d_{0}(x_{n+1})+b(x_n)b_{0}(x_n)\left(b_{1}(x_{n+1})(x_{n+1}-x_n)+f_{0}(x_{n+1})\right)\right] \\
		&&\quad+\left(b_{0}(x_n)b_{0}(x_{n+1})^{2}+2s_{n}\left(l_{0}(x_{n+1})-b_{0}(x_n)\kappa_{0}(x_{n+1})\right)\right)\Bigg]\,h_{2}\Big(\frac{2}{\eps}s_{n}\Big) \\
		&&+\eps^{4} e^{\frac{2 i}{\eps}\phi(x_n)}\Bigg[\frac{1}{2}\left[e_{0}(x_{n+1})+d_{1}(x_{n+1})+b(x_n)b_{0}(x_n)\left(g_{0}(x_{n+1})+f_{1}(x_{n+1})\right)\right] \\
		&&\quad+2\left[b_{0}(x_n)b_{0}(x_{n+1})b_{1}(x_{n+1})+\left(l_{0}(x_{n+1})-b_{0}(x_n)\kappa_{0}(x_{n+1})\right)\right]\Bigg]\,h_{3}\Big(\frac{2}{\eps}s_{n}\Big)\,,
\end{eqnarray*}
and the abbreviations
\begin{eqnarray*}
  && c_{0}(x):=\frac{b(x)^{2}b_{0}(x)}{2\phi'(x)},\: c_{1}(x):=\frac{c_{0}'(x)}{2\phi'(x)},\:d_{0}(x):=\frac{c_{0}(x)}{2\phi'(x)},\: d_{1}(x):=\frac{d_{0}'(x)}{2\phi'(x)},\: e_{0}(x):=\frac{c_{1}(x)}{2\phi'(x)} ,\:\\
  && f_{0}(x):=\frac{b_{0}(x)}{2\phi'(x)},\:
  f_{1}(x):=\frac{f_{0}'(x)}{2\phi'(x)},\: g_{0}(x):=\frac{b_{1}(x)}{2\phi'(x)},\: 
  \kappa_{0}(x):=\frac{b(x)b_{1}(x)}{2\phi'(x)},\: 
%  \kappa(x):=b(x)b_{1}(x),\: \kappa_{0}(x):=\frac{\kappa(x)}{2\phi'(x)},\: 
  l_{0}(x):=\frac{b(x)b_{0}(x)b_{1}(x)}{2\phi'(x)}\,.
%  l(x):=b(x)b_{0}(x)b_{1}(x),\: l_{0}(x):=\frac{l(x)}{2\phi'(x)}\,.
\end{eqnarray*}

\medskip
For these two schemes the following error estimates were proven in \cite{ABN, ADK}:
\begin{theorem}\label{thm:error-est}
Let the coefficient $a\in C^\infty[0,1]$ satisfy $a(x)\ge a_0>0$ in $[0,1]$, and let $0<\eps\le\eps_0$ (for some $0<\eps_0\le1$ such that $\phi'(x)\ne0$ for all $x\in[0,1]$ and $0<\eps\le\eps_0$). Then the global errors of the schemes WKB2 and WKB3 satisfy respectively
%\begin{enumerate}
%\item[(a)]
\begin{equation}\label{error_Z_WKB2}
\|Z_{n}-Z(x_{n})\|_{} \le C \eps^3 h^2\,,\quad \|U_n-U(x_{n})\|_{} \le C {h^{\gamma}
\over \eps} +C \eps^3 h^2\,,\quad n=0,\dots,N\,,
\end{equation} 
%\item[(b)]
\begin{equation}\label{error_Z_WKB3}
\|Z_{n}-Z(x_{n})\|_{} \le C \eps^3 h^3 \max(\eps,h)\,,\quad \|U_n-U(x_n)\|_{} \le C {h^{\gamma}
\over \eps} +C \eps^3 h^3 \max(\eps,h)\,,\quad n=0,\dots,N\,,
\end{equation}
%\end{enumerate}
with $C$ independent of $n$, $h$, and $\eps$.
Here, $\gamma >0$ is the order of the chosen numerical integration method for computing the approximation $\phi_n$ of the phase integral \eqref{phase}, 
\begin{comment}
\begin{equation} \label{PH}
\Phi^{\eps}(x)= \int_0^x \left(\sqrt{a(\tau)} -\eps^2 \beta(\tau) \right) d\tau
\, \left(
\begin{array}{cc}
1&0\\0&-1
\end{array}
\right)
\,.
\end{equation}
\end{comment}
%The discretization step size is $h>0$, 
and $\|.\|$ denotes any vector norm in $\CC^2$.
\end{theorem}
The estimates \eqref{error_Z_WKB2} and \eqref{error_Z_WKB3} include the phase error $|\phi_n-\phi(x_n)|$ only in the backward transformation 
%from $Z_n$ to $U_n$
$U_n=\PP^{-1} e^{{i\over \eps} \PPhi^\eps_n} Z_n\,.$ In \cite{AKU, ADK}, extended error estimates also include the phase error of the analytic transformation from $U$ to $Z$.
For simplicity we used here only a uniform spatial grid; an extension with an adaptive step size controller as well as a coupling to a Runge-Kutta method close to turning points and for the evanescent regime (i.e.\ for $a(x)<0$) is presented in \cite{KAD, ADK}.

%%%%%%%%%%%%%%%%%%%%%%%%%%%%%%%%%%%%%%%%%%%%%%%%%%%%%%%%%%%%%%

\section{NUMERICAL TEST}

\begin{figure}[h]
  \centerline{\includegraphics[width=467pt]{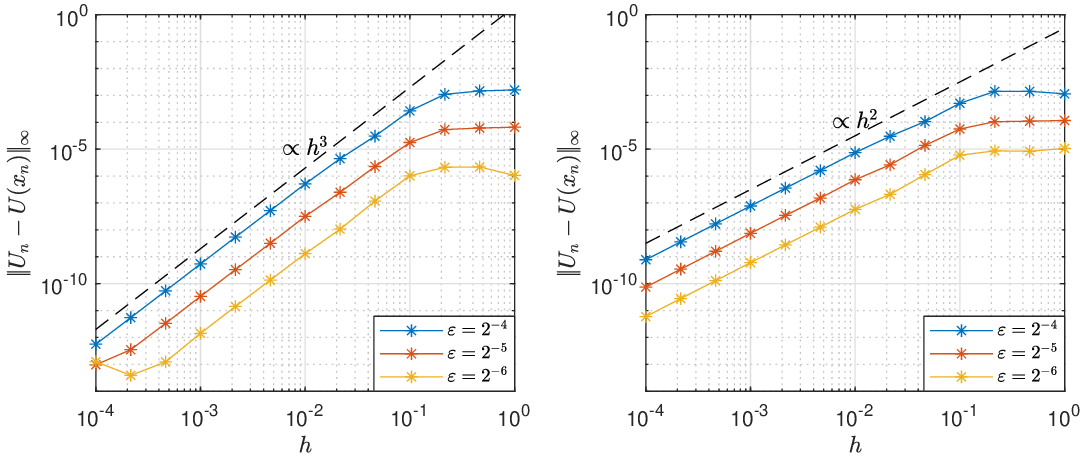}}
  \caption{Log-log plot of the $L^\infty$--error of $U$ as a function of the step size $h$ and for three values of $\eps$, computed with WKB3 (left) and WKB2 (right). The error curve saturates around $10^{-13}$ due to round-off errors.}\label{Fig1}
\end{figure}

We revisit the example from \cite{ABN} with $a(x)=(x+\frac12)^2$. The initial conditions for \eqref{SchIVP} are chosen as $\varphi_0=1$ and $\varphi_1=i$. In Figure \ref{Fig1} we present the $L^\infty$--error of the numerical approximation on $[0,1]$, i.e.\ $\|U_n-U(x_n)\|_\infty$ as a function of the step size $h$ for several values of $\eps$, computed with both WKB3 and WKB2. The error plots are in close agreement with the error estimates \eqref{error_Z_WKB3}, \eqref{error_Z_WKB2}, both when reducing $h$ and $\eps$. Since the phase \eqref{phase} is explicitly computable in this example, the error term $h^\gamma/\eps$ drops out here.

Since the numerical scheme of WKB3 is much more involved than WKB2, and using a lot more function calls, the efficiency gain of WKB3 cannot be inferred only from Figure \ref{Fig1}. But a detailed analysis of the CPU times of both methods at comparable error levels shows a speed-up by up to a factor of 20 for highly accurate computations \cite{ADK}.

% Jannis, 29.7.2022:
% Ich habe als Referenzlösung die Lösung verwendet, welche mir die symbolische Toolbox von Matlab errechnen konnte. Die müsste also "exakt" sein.
%%%%%%%%%%%%%%%%%%%%%%%%%%%%%%%%%%%%%

\begin{comment}

\section{FINAL KEY POINTS TO CONSIDER (FIRST LEVEL HEADING)}
Here are the main points you need to follow (the AIP Publishing author template packages contain comprehensive guidance):

\begin{itemize}
\item Write and prepare your article using the AIP Publishing template.
\item Create a PDF file of your paper (making sure to embed all fonts).
\item Send the following items to your conference organizer:
\begin{itemize}
\item PDF file of your paper
\item Signed Copyright Transfer Agreement
\item (If it applies) Copies of any permissions to re-use copyrighted materials in your article (e.g., figures from books/journals)
\end{itemize}
\end{itemize}

\end{comment}
%%%%%%%%%%%%%%%%%%%%%%%%%%%%%%%%%%%%%%%

% Acknowledgement
\section{ACKNOWLEDGMENTS}
The authors acknowledge support by the projects I3538-N32 and the doctoral school W1245 of the FWF.

% References

\nocite{*}
\bibliographystyle{aipnum-cp}%
%\bibliography{sample}%

\end{document}